\documentclass[11pt,reqno]{amsart}

\usepackage{amsmath,amssymb,amsthm}
\usepackage{graphicx}
\usepackage[a4paper,width=150mm,top=25mm,bottom=25mm,bindingoffset=6mm]{geometry}
\usepackage{fancyhdr}
\allowdisplaybreaks
\usepackage{xcolor}
\usepackage[linktoc=all,colorlinks,linkcolor=teal,citecolor=cyan]{hyperref}
\usepackage{cleveref}
\numberwithin{equation}{section}
\usepackage{comment}
\usepackage{hyperref}

\newcommand{\mb}[1]{\mathbb{#1}}
\newcommand{\mc}[1]{\mathcal{#1}}
\newcommand{\mf}[1]{\mathfrak{#1}}

\renewcommand{\(}{\left(}

\renewcommand{\mod}[1]{\left\lvert#1\right\rvert}

\newtheorem*{thmA}{Theorem A}
\newtheorem*{thmB}{Theorem B}
\newtheorem*{thmC}{Theorem C}

\newtheorem{thm}{Theorem}[section]
\newtheorem{lemma}{Lemma}[section]
\newtheorem{prop}{Proposition}

\newtheorem{remark}{Remark}[section]

\newtheorem*{conj}{Conjecture}

\newtheorem*{proof2}{\emph{Proof of Lemma 3.2}}
\newtheorem{corollary}{Corollary}[section]

\title[Primes in Tuples of Linear Forms in Number Fields and Function Fields]{Primes in Tuples of Linear Forms in Number Fields and Function Fields}
	\author{Habibur Rahaman}
	\address{Department of Mathematics and Statistics, Indian Institute of Science Education and Resaerch Kolkata, Mohanpur, West Bengal - 741246, India.}
	\email{hr21rs044@iiserkol.ac.in, habibr0805@gmail.com}
	
	\thanks{ 2020 \textit{Mathematics Subject Classification.} Primary 11N05, 11N36, 11T06}
 \thanks{\textit{ Key words and phrases.} Distribution of primes, sieve methods, number fields}

\begin{document}
	\begin{abstract}
Following the work of Castillo-Hall-Oliver-Pollack-Thompson who extended Maynard-Tao theorem on admissible tuples to number fields and function fields for tuples with  monic linear forms, here we obtain the Maynard-Tao theorem for  admissible tuples of linear forms with arbitrary leading coefficients in number fields and function fields. Also, we provide some applications of our results.
	\end{abstract}
	\maketitle
	\section{Introduction}


	A set $\mathcal{H}=\{h_1,\ldots,h_k\}$ of distinct integers is called admissible if for any prime $p\in\mathbb{P}$, the set $ \mathcal{H}$ misses at least one residue class $\bmod\  p$, i.e. if $|\mathcal{H}(\bmod\ p)|<p$ for each $p\in\mathbb{P}$. Then we have the following conjecture due to Hardy-Littlewood.
	
	\begin{conj}[Prime k-tuples conjecture]
		For an admissible set $\mathcal{H}=\{h_1,\ldots,h_k\}$ of distinct integers, there are infinitely many $n\in\mb{N}$ such that each of $n+h_i, 1\leq i\leq k$ is a prime.	
	\end{conj}
	 In particular, the above conjecture implies the  twin prime conjecture. Though this conjecture still remains unsolved, recently there has been spectacular progress in this direction. The prime number theorem implies that the average gap between consecutive primes $p_n$ and $p_{n+1}$ is $O(\log p_n)$, where $p_n$ denotes the n-th prime. In \cite{GPY}, Goldston, Pintz and Yildirim proved the following 
	\begin{align*}
		\liminf_{n\to\infty}\frac{p_{n+1}-p_n}{\log p_n}=0,
	\end{align*}
	i.e., the gaps between consecutive primes can be arbitrarily small in comparison with the average prime gap. Zhang \cite{Zhang} extend this result and proved that there are infinitely many consecutive primes with a bounded gap. More precisely, he proved
	\begin{align*}
		\liminf_{n\to\infty}({p_{n+1}-p_n})\leq 70000000.
	\end{align*}
	Maynard \cite{Maynard}, Tao and the Polymath project \cite{Polymath}, managed to reduce the above bound by $246$. Moreover, the method of Maynard and Tao shows that for every $m$, there exist infinitely many tuples of bounded length of $m$ consecutive primes. More precisely, they proved:
	\begin{thmA}[Maynard-Tao]
		Let $m\geq 2$. Then there exists a positive integer $k_0=k_0(m)$ such that for any admissible set $\mathcal{H}=\{h_1,\ldots,h_k\}$ with $k\geq k_0$, there are infinitely many $n$ such that at least $m$ of $n+h_i, 1\leq i\leq k$ are prime.
	\end{thmA}
	\noindent In $2015$, Castillo, Hall, Lemke Oliver, Pollack and Thompson (CHLOPT) \cite{Castillo} extended the Maynard-Tao theorem to number fields and  function fields in the following sense: Let $K$ be a number field with ring of integers $\mathcal{O}_{K}$. We say, the set $\mathcal{H}=\{h_1,\ldots,h_k\}$ of distinct elements in $\mathcal{O}_K$ is admissible if, $|\mathcal{H}(\bmod\ \mathfrak{p})|<|\mathcal{O}_K/\mathfrak{p}|$ for any prime ideal $\mathfrak{p}$ in $\mathcal{O}_K$. Castillo et al. proved the following theorems.
	\begin{thmB}[Theorem 1.1 in \cite{Castillo}]
		Let $m\geq 2$ be a positive integer. Then there exists a positive integer $k_0=k_0(m,K)$ such that, for any admissible set $\mathcal{H}=\{h_1,\ldots,h_k\}$ in $\mathcal{O}_K$ with $k\geq k_0$, there are infinitely many $\alpha\in\mathcal{O}_K$ such that at least $m$ of $\alpha+h_i, 1\leq i\leq k$ are prime.
	\end{thmB}
	
	For the function field setting, we replace $K$ by $\mathbb{F}_q(t)$. Then the ring of integers $\mathcal{O}_K$ is replaced by $\mathbb{F}_q[t]$ where the primes are the monic irreducible polynomials in $\mathbb{F}_q[t]$.
	With this notion, we say $\mathcal{H}=\{h_1,\ldots,h_k\}$ of distinct elements in $\mathbb{F}_q[t]$ is admissible if, for each monic irreducible polynomial $P(t)$, the set $\mathcal{H}(\bmod\ P(t))$ is not equal to $\mathbb{F}_q[t]/P(t)$.
	\begin{thmC}[Theorem 1.3 in \cite{Castillo}]
		Let $m\geq 2$ be a positive integer. Then there exists a positive integer $k_0=k_0(m)$ such that, for any admissible set $\mathcal{H}=\{h_1,\ldots,h_k\}$ in $\mathbb{F}_q[t]$ with $k\geq k_0$, there are infinitely many $f\in\mathbb{F}_q[t]$ such that at least $m$ of $f+h_i, 1\leq i\leq k$ are prime. 
	\end{thmC}
	
	\subsection{Main Results}\label{1.1}
	In this paper, we prove the above results for admissible sets of linear forms with some applications analogous to the results of Banks, Freiberg, and Turnage-Butterbaugh \cite{Banks}. In the number field case, let $\mathcal{H}=\{a_1\alpha+h_1,\ldots,a_k\alpha+h_k\}$, where $\alpha,a_i,h_i\in\mathcal{O}_K$ for $1\leq i\leq k$ and $a_i,h_i$ are coprime for $1\leq i\leq k$.  Let $L_i(\alpha)=a_i\alpha+h_i$ and 
	\begin{align*}
		L(\alpha):=\prod_{i=1}^{k}L_i(\alpha).
	\end{align*}
	We say $\mathcal{H}$ is admissible if $L(\alpha)$ has no fixed prime ideal divisor, that is, if for any prime ideal $\mf{p}$ in $\mathcal{O}_K$, there exists some $\alpha_{\mf{p}}\in\mathcal{O}_K$ such that $\mf{p}\nmid L(\alpha_{\mf{p}})$. Here we prove the following theorem.
	\begin{thm}\label{Thm1}
		Let $m\geq 2$ be a positive integer. Then there exists a positive integer $k_0=k_0(m,K)$ such that, for any admissible set $\mathcal{H}=\{a_1\alpha+h_1,\ldots,a_k\alpha+h_k\}$ of linear forms in $\mathcal{O}_K$ with $k\geq k_0$, there are infinitely many $\alpha\in\mathcal{O}_K$ such that at least $m$ of $a_i\alpha+h_i, 1\leq i\leq k$ are prime.
	\end{thm}
	In particular, for $\mathcal{O}_K=\mathbb{Z}$ with $a_i, h_i$ coprime for $ 1\leq i\leq k$, we have the following corollary (see \cite{Granville}, \cite{Maynard}).
	\begin{corollary}[Theorem 6.4 in \cite{Granville}]
		Let $m\geq 2$. Then there exists a positive integer $k_0=k_0(m)$ such that, for any admissible set $\mathcal{H}=\{a_1n+h_1,\ldots,a_kn+h_k\}$ with $k\geq k_0$, there are infinitely many $n$ such that at least $m$ of $a_in+h_i, 1\leq i\leq k$ are prime.
	\end{corollary}
	From the last result, we also have the following corollary.
	\begin{corollary}\label{main_cor}
	Let $b$ and $h\in\mc{O}_K$ are coprime. Also let $\mathfrak{h}$ denote the ideal generated by $h$ with norm $|\mathfrak{h}|\geq 3$. Then for any positive integer $m\geq 2$, there are infinitely many $r\in\mb{N}$ such that 
	\begin{align*}
		\gamma_{r1}\equiv\gamma_{r2}\equiv\cdots\equiv\gamma_{rm}\equiv b (\bmod\ \mathfrak{h}),
	\end{align*}
where $\gamma_{rj}$ denote primes in $\mc{O}_K$.
	\end{corollary}
	
	In case of a function field $\mathbb{F}_q(t)$, we define an admissible set similarly, letting $f,a_i, h_i\in\mathbb{F}_q[t]$ with $a_i$ monic for $1\leq i\leq k$ and by definition, we consider irreducible polynomials $P(t)\in\mathbb{F}_q[t]$ as primes. We then have the following theorem.
	\begin{thm}\label{Thm2}
		Let $m\geq 2$ be a positive integer. Then there exists a positive integer $k_0=k_0(m)$ such that, for any admissible set $\mathcal{H}=\{a_1f+h_1,\ldots,a_kf+h_k\}$ of linear forms in $\mathbb{F}_q[t]$ with $k\geq k_0$, there are infinitely many $f\in\mathbb{F}_q[t]$ such that at least $m$ of $a_if+h_i, 1\leq i\leq k$ are prime.
	\end{thm}
	
\begin{corollary}\label{main_cor_2}
	Let $b$ and $g\in\mb{F}_q[t]$ are coprime and $g$ be monic. Then for any positive integer $m\geq 2$, there are infinitely many $r\in\mb{N}$ such that 
	\begin{align*}
		g_{r1}\equiv g_{r2}\equiv\cdots\equiv g_{rm}\equiv b (\bmod g),
	\end{align*}
	where $g_{rj}$ are monic irreducible polynomials in $\mb{F}_q[t]$.
\end{corollary}

	Two essential ingredients in the proofs of Theorem~\ref{Thm1} and Theorem~\ref{Thm2} are Proposition~\ref{1st_prop} and Proposition~\ref{2nd_prop}. The proofs of these propositions rely on Lemma~\ref{Key_lemma1} and Lemma~\ref{Key_lemma2}, wherein lie the main ideas which we incorporate in this paper in order to use Maynard's technique and also CHLOPT's technique in the setting of linear forms. 
	
	\section{Preliminaries}\label{sec2}
	In this section, in a similar manner as in \cite{Castillo}, we set up a general  method which will help us to prove both Theorem~\ref{Thm1} and Theorem~\ref{Thm2}         
	simultaneously. 
	\subsection{Notations and Definitions }
	\noindent Let $A$ be a Dedekind domain, which will be $\mathcal{O}_K$ for a number field $K$ and for a function field $\mathbb{F}_q(t)$, we will take $A=\mathbb{F}_q[t]$.\\
	
	\noindent For $A=\mathbb{F}_q[t]$ and for any natural number $N$, let $A(N)$ be the collection of monic polynomials in $\mathbb{F}_q[t]$ with norm $N$. For $A=\mathcal{O}_K$, let $A_0(N)$ be the set of $\alpha\in\mathcal{O}_K$ such that for each real embedding $\sigma:K\to \mathbb{C}, 0<\sigma(\alpha)\leq N$ and $|\sigma(\alpha)|\leq N$ for each complex embedding. Define $A(N):= A_0(2N)\setminus A_0(N)$.
	
	\begin{remark}
		For $K=\mathbb{Q}$, we have $\mathcal{O}_K=\mathbb{Z}$ and $\sigma$ is only the identity map. Therefore $A_0(2N)=[1,2N]$ and $A_0(N)=[1,N]$. So, $A(N)=(N,2N]$, which gives the Maynard-Tao case.
	\end{remark}
	
For any nonzero ideal $\mathfrak{q}\subseteq A$, we define the norm of $\mathfrak{q}$ by $|\mathfrak{q}|:=|A/\mathfrak{q}|$ and the Euler Phi function $\phi(\mathfrak{q}):=|(A/\mathfrak{q})^\times|$, which counts the cardinality of the multiplicative group $(A/\mathfrak{q})^\times$ and  the M\H{o}bius function  $\mu$ by
	\begin{align*}
		\mu(\mathfrak{p}):=
		\begin{cases}
			(-1)^r &\text{ if $\mathfrak{q}=\mathfrak{p}_1.\ldots.\mathfrak{p}_r$, where $\mathfrak{p}_i$ are distinct prime ideals in $A$,}\\
			0 & \text{ otherwise;}
		\end{cases}
	\end{align*}
 which can be defined since $A$ is a Dedekind domain.
	Also, for any $s\in\mathbb{C}$, we define the zeta function of $A$ by
	\begin{align}\label{zeta_function}
		\zeta_A(s):=\sum_{\mathfrak{q}\subseteq A}\frac{1}{|\mathfrak{q}|^{s}},
	\end{align}
	where $|\mathfrak{q}|$ denotes the norm of the ideal $\mathfrak{q}$.
	For $A=\mathcal{O}_K$, this is the usual Dedekind zeta function of $K$ and  for $A=\mathbb{F}_q[t]$, as there are exactly $q^d$ monic polynomials of degree $d$ in $A$, we have 
	\begin{align*}
		\zeta_A(s)=\frac{1}{1-q^{1-s}},
	\end{align*}
	(see page 11 in \cite{Rosenbook}).
	Let $A(N;\mathfrak{q},\alpha_0)$ be the set of $\alpha\in A(N)$ with $\alpha\equiv \alpha_0(\bmod\ \mathfrak{q})$. Then we note that 
	\begin{align*}
		|A(N,\mathfrak{q},\alpha_0)|=\frac{|A(N)|}{|\mathfrak{q}|}+O(|\partial A(N,\mathfrak{q},\alpha_0)|),
	\end{align*}
	where $|A(N)|$ denotes the cardinality of the set $A(N)$ (see page 2844 in \cite{Castillo}) and
	\begin{align}\label{error}
		|\partial A(N;\mathfrak{q},\alpha_0)|\ll 
		\begin{cases}
			1 &\text{ if $A=\mathbb{Z}$  or  $A=\mathbb{F}_q[t]$,}\\
			1+(\frac{|A(N)|}{|\mathfrak{q}|})^{1-\frac{1}{d}}  &\text{ if $A=\mathcal{O}_K$ with $[K,\mathbb{Q}]=d$. }
		\end{cases}
	\end{align}
	
	\noindent Let $P(N)$ be the set of primes in $A(N)$, that is $P(N)=P\cap A(N)$, where $P$ is the set of primes in $A$. Then we have the Prime Number Theorem for $P(N;\mathfrak{q},\alpha_0)$, where $\alpha_0$ and $\mathfrak{q}$ are coprime, as
	\begin{align}\label{Prime_number_thm}
		|P(N;\mathfrak{q},\alpha_0)|=\frac{|P(N)|}{\phi(\mathfrak{q})}+\mathcal{E}(N;\mathfrak{q},\alpha_0),
	\end{align}
	where the error term $\mathcal{E}(N;\mf{q},\alpha_0)$ can be bounded on average by the well known Bombieri-Vinogradov theorem and it's generalization to other fields. To see that, we say the set of primes $P$ in $A$ has level of distribution $\theta>0$ if, for any $B>0$, we have 
	\begin{align*}
		\sum_{|\mathfrak{q}|\leq |A(N)|^{\theta}}\max_{\substack{\alpha_0(\bmod\mathfrak{q})\\(\alpha_0,\mathfrak{q})=1}}|\mathcal{E}(N;\mathfrak{q},\alpha_0)|\ll_B \frac{|A(N)|}{(\log N)^B}.
	\end{align*}
	In particular, if we take $A=\mathbb{Z}$, the Bombieri-Vinogradov theorem says that the primes has level of distribution $\theta$ for any $\theta<1/2$. 
	For a general number field $K$, we have the following result due to Hinz \cite{Hinz_Bombieri}.
	\begin{thm}[Hinz]\label{Hinz}
		Let $K/\mb{Q}$ be a number field with signature $(r_1, r_2)$. If $K$ is totally real, i.e, if $r_2=0$, then the set of primes $P$ in $\mathcal{O}_K$ has level of distribution $\theta$ for any $\theta<1/2$. In general $\theta<\frac{1}{r_2+\frac{5}{2}}$. 
	\end{thm}
	\begin{thm}[Hayes, \cite{Hayes}]\label{Hayes}
		If $A=\mathbb{F}_q[t]$, then the set of primes $P$ has level of distribution $\theta$ for any $\theta<1/2$. Also, we have
		\begin{align*}
			\max_{\substack{\alpha_0(\bmod\mathfrak{q})\\(\alpha_0,\mathfrak{q})=1}}|\mathcal{E}(N;\mathfrak{q},\alpha_0)|\ll (\log2|\mathfrak{q}|)|A(N)|^{1/2}.
		\end{align*}
	\end{thm}
	
	\subsection{Sieve Setting}\label{2.2}
	Let us choose a positive integer $N$ and $D_0=\log\log\log N$. Also, let $\mathfrak{w}:=\displaystyle\prod_{\substack{|\mathfrak{p}|<D_0}}\mathfrak{p}$. Clearly we have
	\begin{align}\label{W_estimate}
		|\mathfrak{w}|\ll (\log\log N)^2.
	\end{align}
	Let $\mathcal{H}=\{a_1\alpha+h_1,\ldots,a_k\alpha+h_k\}$ be an admissible set in $A$. So for any prime ideal $\mathfrak{p}$ in $A$, there is $\alpha_\mathfrak{p}\in A$ such that $\mathfrak{p}\nmid L(\alpha_\mathfrak{p})$. Solving $\alpha\equiv \alpha_\mathfrak{p} (\bmod\ \mathfrak{p})$ for each $\mathfrak{p}\mid\mathfrak{w}$, by Chinese Remainder Theorem there exists some $v_0 (\bmod\ \mathfrak{w})$ such that each $a_iv_0+h_i$ lies in $(A/\mathfrak{w})^{\times}$.\\
	Let $\theta$ be the level of distribution of primes in $A$ and $R=\mod{A(N)}^{\theta/2 -\delta}$ for any fixed small $\delta>0$. Let  $\lambda_{\mathfrak{d}_1,\ldots,\mathfrak{d}_k}$ be suitably chosen weights supported on the tuples $(\mathfrak{d}_1,\ldots,\mathfrak{d}_k)$ of ideals in $A$ with $\mathfrak{d}:=\prod_{i=1}^{k}\mathfrak{d}_i$ squarefree, $|\mathfrak{d}|<R$ and coprime to $\mathfrak{w}$. Next we define two sums, the evaluation of which are crucial for our main results.\\
	
	Let 
	\begin{align}\label{def_S_1}
	\hspace{-4cm}S_1:=\sum_{\substack{\alpha\in A(N)\\\alpha\equiv v_0(\bmod \mathfrak{w})}}\left(\sum_{\substack{\mathfrak{d}_1,\ldots,\mathfrak{d}_k\\a_i\alpha+h_i\equiv 0(\bmod\mathfrak{d}_i)\forall i}}\lambda_{\mathfrak{d}_1,\ldots,\mathfrak{d}_k}\right)^2,
  \end{align}
  \begin{align}\label{def_S_2}
		\hspace{5pt}S_2:=\sum_{\substack{\alpha\in A(N)\\\alpha\equiv v_0(\bmod \mathfrak{w})}}\left(\sum_{i=1}^{k}\chi_P(a_i\alpha+h_i)\right)\left(\sum_{\substack{\mathfrak{d}_1,\ldots,\mathfrak{d}_k\\a_i\alpha+h_i\equiv 0(\bmod\mathfrak{d}_i)\forall i}}\lambda_{\mathfrak{d}_1,\ldots,\mathfrak{d}_k}\right)^2,
	\end{align}
	where $\chi_P(\cdot)$ denotes the characteristic function on the set of primes $P$ in $A$. 
	\begin{remark}
		If we can show that $S_2-\rho S_1>0$ for some $\rho>0$, then there exists at least one $\alpha\in A(N)$ such that at least $\lfloor \rho +1\rfloor$ of $a_1\alpha+h_1,\ldots,a_k\alpha+h_k, 1\leq i\leq k$ are prime.
	\end{remark}
	For evaluating the above sums, we use the following propositions which we prove in the last section.
	
	\begin{prop}\label{1st_prop}
		Let the primes $P$ in $A$ have level of distribution $\theta>0$, and $R=|A(N)|^{\theta/2 -\delta}$ for some fixed $\delta>0$. Let $F:[0,1]^k\to\mathbb{R
		}$ be a piecewise differentiable function supported on $\mathcal{R}_k:=\{(t_1,\ldots,t_k)\in[0,1]^k:\sum_{i=1}^{k} t_i\leq1\}$, and let the sieving weights $\lambda_{\mf{d}_1,\ldots,\mf{d}_k}$ be defined in terms of $F$ as 
		\begin{align*}
			\lambda_{\mathfrak{d}_1,\ldots,\mathfrak{d}_k}=\left(\prod_{i=1}^{k}\mu(\mathfrak{d}_i)\mod{\mathfrak{d}_i}\right)\sum_{\substack{\mathfrak{r}_1,\ldots,\mathfrak{r}_k\\\mathfrak{d}_i|\mathfrak{r}_i\forall i\\(\mathfrak{r}_i,\mathfrak{w})=1\forall i}}\frac{\mu\left(\prod_{i=1}^{k}\mathfrak{r}_i\right)^2}{\prod_{i=1}^{k}\phi(\mathfrak{r}_i)}F\left(\frac{\log \mod{\mathfrak{r}_1}}{\log R},\ldots,\frac{\log \mod{\mathfrak{r}_k}}{\log R}\right),
		\end{align*}
		whenever $\left(\prod_{i=1}^{k}\mathfrak{d}_i,\mathfrak{w}\right)=1$ and $\lambda_{\mathfrak{d}_1,\ldots,\mathfrak{d}_k}=0$ otherwise. Then we have
		\begin{align*}
			&S_1=\frac{(1+o(1))\phi(\mathfrak{w})^k \mod{A(N)}(c_A\log R)^k}{\mod{\mathfrak{w}}^{k+1}}I_k(F),\\
			&S_2= \frac{(1+o(1))\phi(\mathfrak{w})^k \mod{P(N)}(c_A\log R)^{k+1}}{\mod{\mathfrak{w}}^{k+1}}\sum_{m=1}^{k}J^{(m)}_k(F),
		\end{align*}
		where $c_A$ is the residue of $\zeta_A(s)$ at $s=1$, and $I_k(F)$ and $J_k(F)$ are the following integrals
		\begin{align}\label{I_k}
			I_k(F)&:=\int_{0}^{1}\ldots\int_{0}^{1}F(t_1,\ldots,t_k)^2 dt_1\ldots dt_k
		\end{align}
		and
		\begin{align}\label{J_k}
			J^{(m)}_k(F)&:=\int_{0}^{1}\ldots\int_{0}^{1}\left(\int_{0}^{1}F(t_1,\ldots,t_k)dt_m\right)^2 dt_1\ldots dt_{m-1} dt_{m+1}\ldots dt_k.
		\end{align}
	\end{prop}
	\begin{prop}\label{2nd_prop}
		Let $\mathcal{H}=\{a_1\alpha+h_1,\ldots,a_k\alpha+h_k\}$ be admissible set in $A$. Also let $\mc{R}_k,I_k(F)$ and $J^{(m)}_k(F)$ be as in Proposition 1 and  $\mathcal{S}_k$ be set of all piecewise differentiable functions $F:[0,1]^k\to\mathbb{R}$ supported on $\mathcal{R}_k$ with $I_k(F)\neq0$ and $J^{(m)}_k(F)\neq0$ for each $m, 1\leq m\leq k$.
		Let 
		\begin{align*}
			M_k=\sup_{F\in\mathcal{S}_k}\frac{\sum_{m=1}^{k}J^{(m)}_k(F)}{I_k(F)},\ \ \ \ r_k=\left\lceil\frac{\theta M_k}{2}\right\rceil.
		\end{align*}
		Then there are infinitely many integers $\alpha\in A$ such that at least $r_k$ of the $a_i\alpha+h_i, 1\leq i\leq k$ are prime.
	\end{prop}
	For sufficiently large $k$, we have a lower bound of the above $M_k$ due to Maynard \cite{Maynard}.
	\begin{lemma}[Proposition 4.3 in \cite{Maynard}]\label{M_k}
		For sufficiently large $k$, we have
		\begin{align*}
			M_k> \log k-2\log\log k-2.
		\end{align*}
	\end{lemma}
	Also, we note the following lemmas as required to prove the above Propositions.
	\begin{lemma}[Theorem~2.2 in \cite{Rosenbook}]\label{PNT_functionfield}
		Let $a(n)$ denote the number of monic irreducible polynomials in $\mathbb{F}_q[t]$ of degree $n$. Then
		\begin{align*}
			a(n)=\frac{q^n}{n}+O\left(\frac{q^{n/2}}{n}\right).
		\end{align*}
	\end{lemma}
	\begin{lemma}[Mitsui's Generalized Prime Number Theorem, \cite{Mitsui}]\label{GPNT}
		Let $A_0(N)$ be defined as in the beginning of this section for a number field $K/\mathbb{Q}$ with degree $d=r_1+2r_2$. Then 
		\begin{align*}
			\mod{P_0(N)}=m_K\underbrace{\int_{2}^{N}\ldots\int_{2}^{N}}_{r_1\text{ times }}\underbrace{\int_{2}^{N^2}\ldots\int_{2}^{N^2}}_{r_2\text{ times }}\frac{du_1\ldots du_{r_1+r_2}}{\log(u_1\ldots u_{r_1+r_2})}
		+O\left(N^d e^{-c\sqrt{\log( N^d)}}\right),
		\end{align*}
		for some positive real number $c$ and	$m_K:=\frac{w_K}{2^{r_1}h_K R_K}$, where $P_0(N)$ is the set of primes in $A_0(N)$, $w_K$ is the number of roots of unity contained in $K$, $h_K$ is the class number of $K$ and $R_K$ is the regulator of $K$.
	\end{lemma}
	\begin{lemma}[Lemma~6 in \cite{Hinz2}]\label{Value_of_int}
		Let $y_1, y_2,\ldots, y_r$ be real numbers satisfying 
		\begin{align*}
			2\leq y_i\leq y_j^{u},\ \ \ i,j=1,2,\ldots,r,
		\end{align*}
		for a fixed positive real number $u$. Then we have 
		\begin{align*}
			\int_{2}^{y_1}\int_{2}^{y_2}\ldots\int_{2}^{y_r}\frac{dx_1 dx_2\ldots dx_r}{\log(x_1 x_2\ldots x_r)}=\frac{y_1 y_2\ldots y_r}{\log(y_1 y_2\ldots y_r)}+O_u\left(\frac{y_1 y_2\ldots y_r}{\log^2(y_1 y_2\ldots y_r)}\right).
		\end{align*}
	\end{lemma}
	\section{Preparation for the main results}
	In this section, we prove the above propositions. Lets begin with the following lemma, which we get by diagonalizing the quadratic form in $S_1$ (see \eqref{def_S_1}).
	
	\begin{lemma}\label{Key_lemma1}
		For ideals $\mf{r}_1,\ldots,\mf{r}_k$, we set
		\begin{align}\label{def_y}
			y_{\mathfrak{r}_1,\ldots,\mathfrak{r}_k}=\left(\prod_{i=1}^{k}\mu(\mathfrak{r}_i)\phi(\mathfrak{r}_i)\right)\sum_{\substack{\mathfrak{d}_1,\ldots,\mathfrak{d}_k\\\mathfrak{r}_i\mid \mathfrak{d}_i,\forall i}}\frac{\lambda_{\mathfrak{d}_1,\ldots,\mathfrak{d}_k}}{\prod_{i=1}^{k}\mod{\mathfrak{d}_i}},
		\end{align}
		and $y_{\max}=\displaystyle\sup_{\mathfrak{r}_1,\ldots,\mathfrak{r}_k}|y_{\mathfrak{r}_1,\ldots,\mathfrak{r}_k}|$. Then
		\begin{align}
			S_1=\frac{\mod{A(N)}}{\mod{\mathfrak{w}}}\sum_{\mathfrak{r}_1,\ldots,\mathfrak{r}_k}\frac{y^2_{\mathfrak{r}_1,\ldots,\mathfrak{r}_k}}{\prod_{i=1}^{k}\phi(\mathfrak{r}_i)}+O\left(\frac{y^2_{\max}\phi(\mathfrak{w})^k \mod{A(N)}(\log R)^k}{\mod{\mathfrak{w}}^{k+1}D_0}\right).
		\end{align}
	\end{lemma}
	
	\begin{remark}
		The above change of variables is invertible (see page 393 in \cite{Maynard}). For squarefree $\mathfrak{d}=\prod_{i=1}^{k}\mathfrak{d}_i$, we have
		\begin{align*}
			\sum_{\substack{\mathfrak{r}_1,\ldots,\mathfrak{r}_k\\\mathfrak{d}_i\mid \mathfrak{r}_i,\forall i}}&\frac{y_{\mathfrak{r}_1,\ldots,\mathfrak{r}_k}}{\prod_{i=1}^{k}\phi({\mathfrak{r}_i})}=\sum_{\substack{\mathfrak{r}_1,\ldots,\mathfrak{r}_k\\\mathfrak{d}_i\mid \mathfrak{r}_i,\forall i}}\left(\prod_{i=1}^{k}\mu(\mathfrak{r}_i)\right)\sum_{\substack{\mathfrak{e}_1,\ldots,\mathfrak{e}_k\\\mathfrak{r}_i\mid \mathfrak{e}_i,\forall i}}\frac{\lambda_{\mathfrak{e}_1,\ldots,\mathfrak{e}_k}}{\prod_{i=1}^{k}\mod{\mathfrak{e}_i}}\\
			&=\sum_{\mathfrak{e}_1,\ldots,\mathfrak{e}_k}\frac{\lambda_{\mathfrak{e}_1,\ldots,\mathfrak{e}_k}}{\prod_{i=1}^{k}\mod{\mathfrak{e}_i}}\sum_{\substack{\mathfrak{r}_1,\ldots,\mathfrak{r}_k\\\mathfrak{d}_i\mid \mathfrak{r}_i,\forall i\\\mathfrak{r}_i\mid \mathfrak{e}_i,\forall i\\}}\prod_{i=1}^{k}\mu(\mathfrak{r}_i)
			=\frac{\lambda_{\mathfrak{d}_1,\ldots,\mathfrak{d}_k}}{\prod_{i=1}^{k}\mu_i(\mathfrak{d}_i)\mod{\mathfrak{d}_i}}.
		\end{align*}
	Choosing suitable $y_{\mathfrak{r}_1,\ldots,\mathfrak{r}_k}$ in terms of the function $F$, we get the sieve weights given in Proposition 1. Also, we note that $\lambda_{\max}\ll y_{\max}(\log R)^k$.
	\end{remark}
	\noindent\emph{Proof of Lemma~\ref{Key_lemma1}.}
	Let us recall from \eqref{def_S_1} that 
	\begin{align*}
		S_1=\sum_{\substack{\alpha\in A(N)\\\alpha\equiv v_0(\bmod \mathfrak{w})}}\left(\sum_{\substack{\mathfrak{d}_1,\ldots,\mathfrak{d}_k\\a_i\alpha+h_i\equiv 0(\bmod\mathfrak{d}_i)\forall i}}\lambda_{\mathfrak{d}_1,\ldots,\mathfrak{d}_k}\right)^2.
	\end{align*}
	Since the sum is over a finite set, we expand the inner sum and  then interchanging the order of summation, we get 
	\begin{align*}
		S_1=\sum_{\substack{\mathfrak{d}_1,\ldots,\mathfrak{d}_k\\\mathfrak{e}_1,\ldots,\mathfrak{e}_k}}\lambda_{\mathfrak{d}_1,\ldots,\mathfrak{d}_k}\lambda_{\mathfrak{e}_1,\ldots,\mathfrak{e}_k}\sum_{\substack{\alpha\in A(N)\\\alpha\equiv v_0(\bmod \mathfrak{w})\\a_i\alpha+h_i\equiv 0(\bmod [\mathfrak{d}_i,\mathfrak{e}_i] )\forall i}}1
	\end{align*}
	Now $a_i\alpha+h_i\equiv 0(\bmod\ [\mathfrak{d}_i,\mathfrak{e}_i])$ is equivalent to $\alpha\equiv-a^{-1}_ih_i(\bmod\ [\mathfrak{d}_i,\mathfrak{e}_i])$,  as $a_i, h_i$ are coprime. Thus by Chinese Remainder Theorem, we can write the inner sum in a single residue class $\alpha_0$ modulo $\mathfrak{q}:=\mathfrak{w}\prod_{i=1}^{k}[\mathfrak{d}_i,\mathfrak{e}_i]$, provided $\mathfrak{w}, [\mathfrak{d}_1,\mathfrak{e}_1],\ldots,[\mathfrak{d}_k,\mathfrak{e}_k]$ are pairwise coprime, otherwise the inner sum is zero. Also, for the coprime case, by \eqref{error} the inner sum is equal to 
	\begin{align*}
		\frac{|A(N)|}{\mod{\mathfrak{q}}}+O(\mod{\partial A(N;\mathfrak{q},\alpha_0)})=\frac{|A(N)|}{\mod{\mathfrak{q}}}
		+O\left(\frac{\mod{A(N)}}{\mod{\mathfrak{q}}}\right)^{1-l},
	\end{align*}
	for some $l>0$, provided $\mod{A(N)}\geq \mod{\mathfrak{q}}$, which we can get by taking $N$ sufficiently large. Therefore 
	\begin{align}\label{S_1_beffore_manipulation}
		S_1=\frac{\mod{A(N)}}{\mod{\mathfrak{w}}}\sideset{}{'}\sum_{\substack{\mathfrak{d}_1,\ldots,\mathfrak{d}_k\\\mathfrak{e}_1,\ldots,\mathfrak{e}_k}}\frac{\lambda_{\mathfrak{d}_1,\ldots,\mathfrak{d}_k}\lambda_{\mathfrak{e}_1,\ldots,\mathfrak{e}_k}}{\prod_{i=1}^{k}\mod{[\mathfrak{d}_i,\mathfrak{e}_i]}}+O\left(\mod{A(N)}^{1-l}\sideset{}{'}\sum_{\substack{\mathfrak{d}_1,\ldots,\mathfrak{d}_k\\\mathfrak{e}_1,\ldots,\mathfrak{e}_k}}\frac{|\lambda_{\mathfrak{d}_1,\ldots,\mathfrak{d}_k}\lambda_{\mathfrak{e}_1,\ldots,\mathfrak{e}_k}|}{\mod{\mathfrak{q}}^{1-l}}\right),
	\end{align}
	where $\sideset{}{'}\sum$ denotes that the sum is over squarefree $\mathfrak{q}$. We note that the above error term is
	\begin{align*}
		\ll\lambda_{\max}^2\mod{A(N)}^{1-l}\sideset{}{'}\sum_{\mod{\mathfrak{q}}<\mod{\mathfrak{w}}R^2}\frac{\mu^2(\mathfrak{q})\tau_{3k}(\mathfrak{q})}{\mod{\mathfrak{q}}^{1-l}},
	\end{align*}
	since for each $\mathfrak{q}$, there are at most $\tau_{3k}(\mathfrak{q})$ choices of $\mathfrak{d}_1,\ldots,\mathfrak{d}_k,\mathfrak{e}_1,\ldots,\mathfrak{e}_k$ such that $\mathfrak{q}=\mathfrak{w}\prod_{i=1}^{k}[\mathfrak{d}_i,\mathfrak{e}_i]$, where $\tau_k(\mathfrak{q})$ is the number of ways that $\mathfrak{q}$ can be written as product of $k$ ideals. Thus the error is  
	\begin{align}\label{extra_error}
		&\hspace{-6pt}\ll\hspace{-2pt}\lambda_{\max}^2\mod{A(N)}^{1-l}\hspace{-2pt}(\mod{\mathfrak{w}}R^2)^l\hspace{-6pt}\prod_{\mod{\mathfrak{p}}<\mod{\mathfrak{w}}R^2}\hspace{-5pt}\left(1+\frac{3k}{\mod{\mathfrak{p}}}\right)
		\hspace{-2pt}\ll\hspace{-2pt}\lambda_{\max}^2\mod{A(N)}^{1-l}\hspace{-2pt}(\mod{\mathfrak{w}}R^2)^l(\log(\mod{\mathfrak{w}}R^2))^{3k}\nonumber\\
		& \hspace{-2pt}\ll\lambda_{\max}^2\mod{A(N)}(\log R)^{4k}\left(\frac{R^2}{\mod{A(N)}}\right)^l
		\hspace{-2pt}\ll y^2_{\max}\mod{A(N)}(\log R)^{6k}\frac{1}{|A(N)|^{(1-\theta +2\delta)l}},
	\end{align}
	 where in the second line is due to Mertens' theorem for global fields (see \cite{Rosen2}), third and forth lines are due to \eqref{W_estimate} and $\lambda_{\max}\ll y^2_{\max}(\log R)^k$, respectively. The above error term is dominated by the error stated in this lemma as $\theta<1$ and $l>0$.\\
  Now by following Maynard's (see page 394 in \cite{Maynard}),  using $\frac{1}{|[\mathfrak{d}_i,\mathfrak{e}_i]|}=\frac{1}{|\mathfrak{d}_i||\mathfrak{e}_i|}\sum_{\mathfrak{u}_i\mid \mathfrak{d}_i,\mathfrak{e}_i}\phi(\mathfrak{u}_i)$ and introducing auxiliary variables $\mathfrak{s}_{i,j}$ to uncouple $(\mathfrak{u}_i,\mathfrak{e}_i)=1$, and then using the change of variables \eqref{def_y}, we obtain the equality of the main term in \eqref{S_1_beffore_manipulation} with the following
	\begin{align}\label{Manipulating-sum}
		\frac{\mod{A(N)}}{\mod{\mathfrak{w}}}\sum_{\mathfrak{u}_1,\ldots,\mathfrak{u}_k}\left(\prod_{i=1}^{k}\frac{\mu(\mathfrak{u}_i)^2}{\phi(\mathfrak{u}_i)}\right)\sideset{}{^*}\sum_{\mathfrak{s}_{1,2},\ldots,\mathfrak{s}_{k,k-1}}\left(\prod_{\substack{1\leq i,j\leq k\\ i\neq j}}\frac{\mu(\mathfrak{s}_{i,j})}{\phi(\mathfrak{s}_{i,j})^2}\right)y_{\mathfrak{b}_1,\ldots,\mathfrak{b}_k}y_{\mathfrak{c}_1,\ldots,\mathfrak{c}_k},
	\end{align}
	where $\mathfrak{b}_i:=\mathfrak{u}_i\prod_{j\neq i}\mathfrak{s}_{i,j}, \mathfrak{c}_j:=\mathfrak{u}_j\prod_{i\neq j}\mathfrak{s}_{i,j}$, and $\sideset{}{^*}\sum$ denotes the sum is restricted to $\mathfrak{s}_{i,j}$ such that each $\mathfrak{s}_{i,j}$ is coprime to any other term of $\mathfrak{b}_i$ and $\mathfrak{c}_j$. From the support of $\lambda's$, we can see that $y's$ are also supported on the tuples $(\mathfrak{d}_1,\ldots,\mathfrak{d}_k)$ of ideals in $A$ with $\mathfrak{d}:=\prod_{i=1}^{k}\mathfrak{d}_i$ squarefree, $|\mathfrak{d}|<R$ and coprime to $\mathfrak{w}$. Therefore there is no contribution of the inner sum for $(\mathfrak{s}_{i,j},\mathfrak{w})\neq 1$, and we can take $(\mathfrak{s}_{i,j},\mathfrak{w})= 1$. Then we have only two cases, either $\mathfrak{s}_{i,j}=1, \forall i\neq j$, or $\mod{\mathfrak{s}_{i,j}}>D_0$ for some $i\neq j$. For the later case the contribution of the above main term is 
	\begin{align*}
		&\ll \frac{y^2_{\max}\mod{A(N)}}{\mod{\mathfrak{w}}}\left(\sum_{\substack{\mod{\mathfrak{u}}<R\\(\mathfrak{u},\mathfrak{w})=1}}\frac{\mu(\mathfrak{u})^2}{\phi(\mathfrak{u})}\right)^k \left(\sum_{\mod{\mathfrak{s}_{i,j}}>D_0}\frac{\mu(\mathfrak{s}_{i,j})^2}{\phi(\mathfrak{s}_{i,j})^2}\right)\left(\sum_{\mathfrak{s}\subseteq A}\frac{\mu(\mathfrak{s})^2}{\phi(\mathfrak{s})^2}\right)^{k^2 -k-1}\\
		&\ll \frac{y^2_{\max}\mod{A(N)}\phi(\mathfrak{w})^k (\log R)^k}{\mod{\mathfrak{w}}^{k+1}D_0}.
	\end{align*}
	which is the error stated in the lemma.\\
	\noindent We are left with the case that all $\mathfrak{s}_{i,j}=1$ for all $i\neq j$ in the main term. In this case, the main term is equal to
	\begin{align*}
		\frac{\mod{A(N)}}{\mod{\mathfrak{w}}}\sum_{\mathfrak{r}_1,\ldots,\mathfrak{r}_k}\frac{y^2_{\mathfrak{r}_1,\ldots,\mathfrak{r}_k}}{\prod_{i=1}^{k}\phi(\mathfrak{r}_i)}
	\end{align*}
	and which completes the proof.\qed\\
	
	We now want to evaluate a similar sum for $S_2$. For that, we write 
	$S_2=\sum_{m=1}^{k}S^{(m)}_2$, where
	\begin{align*}
		S^{(m)}_2=\sum_{\substack{\alpha\in A(N)\\\alpha\equiv v_0(\bmod \mathfrak{w})}}\chi_P(a_m\alpha+h_m)\left(\sum_{\substack{\mathfrak{d}_1,\ldots,\mathfrak{d}_k\\(a_i\alpha+h_i)\equiv 0(\bmod\mathfrak{d}_i)\forall i}}\lambda_{\mathfrak{d}_1,\ldots,\mathfrak{d}_k}\right)^2,
	\end{align*}
	and we evaluate $S^{(m)}_2$ in the following lemma.
	\begin{lemma}\label{Key_lemma2}
		Let
		\begin{align}\label{def_ym}
			y^{(m)}_{\mathfrak{r}_1,\ldots,\mathfrak{r}_k}=\left(\prod_{i=1}^{k}\mu(\mathfrak{r}_i)g(\mathfrak{r}_i)\right)\sum_{\substack{\mathfrak{d}_1,\ldots,\mathfrak{d}_k\\\mathfrak{r}_i\mid \mathfrak{d}_i,\forall i\\\mathfrak{d}_m=1}}\frac{\lambda_{\mathfrak{d}_1,\ldots,\mathfrak{d}_k}}{\prod_{i=1}^{k}\phi(\mathfrak{d}_i)},
		\end{align}
		where $g$ is a totally multiplicative function defined on ideals of $A$ as $g(\mathfrak{p})=\mod{\mathfrak{p}}-2$ for each prime ideal $\mathfrak{p}$ of $A$. Also, let $y^{(m)}_{\max}=\displaystyle\sup_{\mathfrak{r}_1,\ldots,\mathfrak{r}_k}|y^{(m)}_{\mathfrak{r}_1,\ldots,\mathfrak{r}_k}|$. Then for any fixed $B>0$, we have
		\begin{align}\label{S_2m}
			S^{(m)}_2=\frac{\mod{P(N)}}{\phi({\mathfrak{w}})}\sum_{\mathfrak{r}_1,\ldots,\mathfrak{r}_k}\frac{(y^{(m)}_{\mathfrak{r}_1,\ldots,\mathfrak{r}_k})^2}{\prod_{i=1}^{k}g(\mathfrak{r}_i)}&+O\left(\frac{(y^{(m)}_{\max})^2\phi(\mathfrak{w})^{k-2} \mod{A(N)}(\log N)^{k-2}}{\mod{\mathfrak{w}}^{k-1}D_0}\right)\nonumber\\&+O\left(\frac{y^2_{\max}\mod{A(N)}}{(\log N)^B}\right).
		\end{align}
	\end{lemma}
\noindent  To ensure that the main term in Lemma~\ref{Key_lemma2} is nonzero, we chose $N$ such that for each prime ideal $\mathfrak{p}$ (in $\mathcal{O}_K$ for Theorem~\ref{Thm1} and in $\mathbb{F}_q[t]$ for Theorem~\ref{Thm2}) with norm larger than $\log\log\log N$,  $\mathfrak{p}\nmid (a_ih_j-a_jh_i)$ for all $i\neq j, 1\leq i,j\leq k$ (see Subsection~\ref{1.1} and Subsection~\ref{2.2}). Then we have reduced the main term in Lemma~\ref{Key_lemma2} by using Mitsui's generalized prime number theorem (see Theorem~\ref{GPNT}). Also to bound the error term in Lemma~\ref{Key_lemma1} and Lemma~\ref{Key_lemma2}, we used Hinz's version (see Theorem~\ref{Hinz}) and Hayes's version (see Theorem~\ref{Hayes}) of Bombieri-Vinogradov theorem. The error term in our case depends on the norms of the principal ideals $\mf{a}_i$ generated by $a_i$, for all $i$ (see Subsection~\ref{1.1}). Here, a priori, it may seem that the error term would exceed the main term. However, we managed to control the error by taking $N$ sufficiently large (see Subsection~\ref{2.2}) compared to the norms of these ideals.
	\begin{proof2} \textnormal{Expanding the inner sum and interchanging the order of summation, we have}
		\begin{align*}
			S^{(m)}_2=\sum_{\substack{\mathfrak{d}_1,\ldots,\mathfrak{d}_k\\\mathfrak{e}_1,\ldots,\mathfrak{e}_k}}\lambda_{\mathfrak{d}_1,\ldots,\mathfrak{d}_k}\lambda_{\mathfrak{e}_1,\ldots,\mathfrak{e}_k}\sum_{\substack{\alpha\in A(N)\\\alpha\equiv v_0(\bmod \mathfrak{w})\\a_i\alpha+h_i\equiv 0(\bmod [\mathfrak{d_i},\mathfrak{e_i}]),\forall i}}\chi_P(a_m\alpha+h_m).
		\end{align*}
		\textnormal{As in Lemma~\ref{Key_lemma1}, the inner sum can be written as a single residue class $\alpha_0$ coprime to $\mathfrak{q}:=\mathfrak{w}\prod_{i=1}^{k}[\mathfrak{d}_i,\mathfrak{e}_i]$, provided $\mathfrak{q}$ is squarefree. As $[\mathfrak{d}_m,\mathfrak{e}_m]\mid (a_m\alpha+h_m)$, if $|[\mathfrak{d}_m,\mathfrak{e}_m]|>1$, for sufficiently large $N$ the above inner sum is zero. Therefore taking $N$ large, we can write the above sum as }
		\begin{align*}
			S^{(m)}_2=\sideset{}{'}\sum_{\substack{\mathfrak{d}_1,\ldots,\mathfrak{d}_k\\\mathfrak{e}_1,\ldots,\mathfrak{e}_k\\
					\mathfrak{d}_m=\mathfrak{e}_m=1}}\lambda_{\mathfrak{d}_1,\ldots,\mathfrak{d}_k}\lambda_{\mathfrak{e}_1,\ldots,\mathfrak{e}_k}\sum_{\substack{\alpha\in A(N)\\\alpha\equiv \alpha_0(\bmod \mathfrak{q})}}\chi_P(a_m\alpha+h_m),
		\end{align*}
  \textnormal{where $\sideset{}{'}\sum$ denotes that the sum is over squarefree $\mf{q}$.\\}
	\textnormal{ We claim that $a_m\alpha_0+h_m$ and $\mathfrak{a}_m\mathfrak{q}$ are coprime, where $\mf{a}_m$ is the principal ideal generated by $a_m$. Suppose they are not coprime and let $\mathfrak{p}\mid (a_m\alpha_0+h_m,\mathfrak{a}_m\mathfrak{q})$. As $\mf{p}\mid(a_m\alpha_0+h_m)$, if $\mathfrak{p}|\mathfrak{a}_m$ then $\mathfrak{p}|h_m$, which contradicts that $a_m$ and $h_m$ are coprime. Therefore $\mf{p}\nmid\mf{a}_m$. Also $\mathfrak{p}\nmid \mathfrak{w}$ as $a_m\alpha_0+h_m\equiv a_mv_0+h_m(\bmod\ \mathfrak{w})$ with $a_mv_0+h_m$ and $\mathfrak{w}$ are coprime. So $\mathfrak{p}\mid (a_m\alpha_0+h_m,\prod_{i=1}^{k}[\mathfrak{d}_i,\mathfrak{e}_i])$. Now, if $\mathfrak{p}|[\mathfrak{d}_i,\mathfrak{e}_i]$ for some $i\neq m$, then $a_m\alpha_0+h_m\equiv0(\bmod\ \mathfrak{p})$ and $a_i\alpha_0+h_i\equiv0(\bmod\ \mathfrak{p})$. As $\mathfrak{p}\nmid \mathfrak{w}$, $|\mathfrak{p}|>D_0$. Therefore, for sufficiently large $N$, we have $\alpha_0\equiv -a^{-1}_mh_m(\bmod\ \mathfrak{p})$ and $\alpha_0\equiv -a^{-1}_ih_i(\bmod\ \mathfrak{p})$, and combining this, we have  $a^{-1}_mh_m\equiv a^{-1}_ih_i(\bmod\ \mathfrak{p})$, that is $a_ih_m- a_mh_i\equiv 0(\bmod\ \mathfrak{p})$. Thus $\mathfrak{p}\mid (a_ih_m- a_mh_i)$, which leads to a contradiction for large $N$. This proves the claim and therefore, the inner sum counts the number of primes in a residue class coprime to the modulus $\mathfrak{a}_m\mathfrak{q}$.}
		
\noindent \textnormal{For a number field $K$, let
		\begin{align}
			A_0(a_m,N):=\{\beta\in &A:   0<\sigma(\beta)\leq |\sigma(a_m)|N,\text{ for each real embedding of $K$ and}\nonumber\\&1\leq\mod{\sigma(\beta)}\leq \mod{\sigma(a_m)}N, \text{ for each complex embedding of $K$}\}\label{A_0(am,N)}
		\end{align}
		and let $A(a_m,N):=A_0(a_m,2N)\setminus A_0(a_m,N)$.}\\
		\textnormal{For a function field $\mathbb{F}_q(t)$, let
		\begin{align*}
			A(a_m,N):=\{\beta\in \mathbb{F}_q[t]: \beta \text { monic and } \mod{\beta}=\mod{\mathfrak{a}_m}N \}.
		\end{align*}
		With the above notations, by \eqref{error} and \eqref{Prime_number_thm} we have for some $l>0$,
		\begin{align}\label{critial_pt}
			\sum_{\substack{\alpha\in A(N)\\\alpha\equiv \alpha_0(\bmod \mathfrak{q})}}\chi_P(a_m\alpha+h_m)=\frac{\mod{P(a_m,N)}}{\phi(\mathfrak{a}_m\mathfrak{q})}&+O\left(\left(\frac{\mod{A(a_m,N)}}{\mod{\mathfrak{a}_m\mathfrak{q}}}\right)^{1-l}\right)\nonumber\\&+O(\mathcal{E}(\mod{\mathfrak{a}_m}N;\mathfrak{a}_m\mathfrak{q},\alpha_1)),
		\end{align}
		where $P(a_m,N)=P\cap A(a_m,N)$. (The second term on the right hand side is necessary for the number field case, for example one can see Lemma~11 in \cite{Anirban}.) Thus
		\begin{align}\label{S_2_before_manipulation}
			S^{(m)}_2\hspace{-4pt}=\hspace{-2pt}\frac{\mod{P(a_m,N)}}{\mod{\phi(\mathfrak{a}_m\mathfrak{w})}}&\hspace{-5pt}\sideset{}{'}\sum_{\substack{\mathfrak{d}_1,\ldots,\mathfrak{d}_k\\\mathfrak{e}_1,\ldots,\mathfrak{e}_k\\
\mathfrak{d}_m=\mathfrak{e}_m=1}}\hspace{-5pt}\frac{\lambda_{\mathfrak{d}_1,\ldots,\mathfrak{d}_k}\lambda_{\mathfrak{e}_1,\ldots,\mathfrak{e}_k}}{\prod_{i=1}^{k}\phi([\mathfrak{d}_i,\mathfrak{e}_i])}\hspace{-2pt}+\hspace{-2pt}O\hspace{-3pt}\left(\sideset{}{'}\sum_{\substack{\mathfrak{d}_1,\ldots,\mathfrak{d}_k\\\mathfrak{e}_1,\ldots,\mathfrak{e}_k}}\mod{\lambda_{\mathfrak{d}_1,\ldots,\mathfrak{d}_k}\lambda_{\mathfrak{e}_1,\ldots,\mathfrak{e}_k}\mathcal{E}(\mod{\mathfrak{a}_m}N;\mathfrak{a}_m\mathfrak{q},\alpha_1)}\hspace{-3pt}\right)\nonumber\\
&+O\left(\lambda^2_{\max}\frac{\mod{A(a_m,N)}^{1-l}}{\mod{\mathfrak{a}_m\mathfrak{w}}^{1-l}}\sideset{}{'}\sum_{\substack{\mathfrak{d}_1,\ldots,\mathfrak{d}_k\\\mathfrak{e}_1,\ldots,\mathfrak{e}_k}}\frac{1}{\prod_{i=1}^{k}\mod{[\mathfrak{d}_i,\mathfrak{e}_i]}^{1-l}}\right).
		\end{align}
		As $a_m$ is fixed in terms of $N$, taking sufficiently large $N$, we have
		\begin{align*}
			\frac{\mod{A(\mod{\mathfrak{a}_m}N)}^{1-l}}{\mod{\mathfrak{a}_m\mathfrak{w}}^{1-l}}\ll \frac{\mod{A(N)}^{1-l}}{\mod{\mathfrak{w}}^{1-l}}
		\end{align*} 
		and by \eqref{extra_error}, second big oh-term is 
		\begin{align*}
			\ll y^2_{\max}\mod{A(N)}(\log R)^{6k}\frac{1}{|A(N)|^{(1-\theta +2\delta)l}},
		\end{align*}
		which is dominated by the second error term stated in the lemma.\\}
		\textnormal{For the first error term, let 
		\begin{align*}
			\mathcal{E}(\mod{\mathfrak{a}_m}N;\mathfrak{a}_m\mathfrak{q})=\max_{\substack{\alpha_1(\bmod\mathfrak{a}_m\mathfrak{q})\\(\alpha_1,\mathfrak{a}_m\mathfrak{q})=1}}\mod{\mathcal{E}(\mod{\mathfrak{a}_m}N;\mathfrak{a}_m\mathfrak{q},\alpha_1)}.
		\end{align*}
		Also, notice that for any $\mathfrak{q}$, there are at most $\tau_{3k}(\mathfrak{q})$ choices of $\mathfrak{d}_1,\ldots,\mathfrak{d}_k,\mathfrak{e}_1,\ldots,\mathfrak{e}_k$ such that $\mathfrak{q}:=\mathfrak{w}\prod_{i=1}^{k}[\mathfrak{d}_i,\mathfrak{e}_i]$,  where $\tau_k(\mathfrak{q})$ is the number of ways that $\mathfrak{q}$ can be written as product of $k$ ideals. So the first error is bounded by 
		\begin{align*}
			\lambda^2_{\max} \sum_{\mod{\mathfrak{r}}<\mod{\mathfrak{a}_m}\mod{\mathfrak{w}}R^2}\mu(\mathfrak{r})^2\tau_{3k}(\mathfrak{r})\mathcal{E}(\mod{\mathfrak{a}_m}N;\mathfrak{r}).
		\end{align*}
		Now $\mod{\mathfrak{r}}<\mod{\mathfrak{a}_m}\mod{\mathfrak{w}}R^2\ll(\log\log N)^2 \mod{A(N)}^{\theta-2\delta}\ll\mod{A(N)}^{\theta-\delta'}$ for some $\delta'>0$, so we can use Theorem~\ref{Hinz} and Theorem~\ref{Hayes}. Then using Cauchy-Schwarz inequality, the trivial bound $\mathcal{E}(\mod{\mathfrak{a}_m}N;\mathfrak{r})\ll \frac{\mod{A(\mod{\mathfrak{a}_m}N)}}{\phi(\mathfrak{a}_m\mathfrak{q})}$, Theorem~\ref{Hinz} and Theorem~\ref{Hayes}, the above error is 
		\begin{align*}
			&\ll \lambda^2_{\max}\left(\sum_{\mod{\mathfrak{r}}<\mod{\mathfrak{a}_m}\mod{\mathfrak{w}}R^2}\mu(\mathfrak{r})^2\tau^2_{3k}(\mathfrak{r})\mathcal{E}(\mod{\mathfrak{a}_m}N;\mathfrak{r})\right)^{1/2}\left(\sum_{\mod{\mathfrak{r}}<\mod{\mathfrak{a}_m}\mod{\mathfrak{w}}R^2}\mu(\mathfrak{r})^2\mathcal{E}(\mod{\mathfrak{a}_m}N;\mathfrak{r})\right)^{1/2}\\
			&\ll \lambda^2_{\max}\left(\sum_{\mod{\mathfrak{r}}<\mod{\mathfrak{a}_m}\mod{\mathfrak{w}}R^2}\mu(\mathfrak{r})^2\tau^2_{3k}(\mathfrak{r})\frac{\mod{A(\mod{\mathfrak{a}_m}N)}}{\phi(\mathfrak{a}_m\mathfrak{q})}\right)^{1/2}\left(\sum_{\mod{\mathfrak{r}}<\mod{\mathfrak{a}_m}\mod{\mathfrak{w}}R^2}\hspace{-6pt}\mu(\mathfrak{r})^2\mathcal{E}(\mod{\mathfrak{a}_m}N;\mathfrak{r})\right)^{1/2}\\
			& \ll y^2_{\max}(\log R)^{2k}\hspace{-4pt}\left(\sum_{\mod{\mathfrak{r}}<\mod{\mathfrak{a}_m}\mod{\mathfrak{w}}R^2}\hspace{-5pt}\mu(\mathfrak{r})^2\tau^2_{3k}(\mathfrak{r})\frac{\mod{A(N)}}{\phi(\mathfrak{q})}\right)^{1/2}\hspace{-7pt}\left(\sum_{\mod{\mathfrak{r}}<\mod{\mathfrak{a}_m}\mod{\mathfrak{w}}R^2}\hspace{-6pt}\mu(\mathfrak{r})^2\mathcal{E}(\mod{\mathfrak{a}_m}N;\mathfrak{r})\right)^{1/2}\\
			&\ll \frac{y^2_{\max}\mod{A(N)}}{(\log N)^B},
		\end{align*}
		for any large $B$ and so, we are done with the error term.\\}
		\textnormal{For the main term on the right hand side of \eqref{S_2_before_manipulation}, let us recall the multiplicative function $g$ stated in the lemma, defined by $g(\mathfrak{p})=\mod{\mathfrak{p}}-2$ and notice that, for squarefree $\mathfrak{d}_i, \mathfrak{e}_i$,
		\begin{align*}
			\frac{1}{\phi([\mathfrak{d}_i,\mathfrak{e}_i])}=\frac{1}{\phi(\mathfrak{d}_i)\phi(\mathfrak{e}_i)}\sum_{\mathfrak{u}_i\mid \mathfrak{d}_i,\mathfrak{e}_i}g(\mathfrak{u}_i).
		\end{align*}
		As $(\mathfrak{d}_i,\mathfrak{e}_i)=1$, we have $\sum_{\mathfrak{s}_{i,j}\mid \mathfrak{d}_i,\mathfrak{e}_i}\mu(\mathfrak{s}_{i,j})=1$. Thus, instead of $(\mathfrak{d}_i,\mathfrak{e}_i)=1$, substituting the expression $\sum_{\mathfrak{s}_{i,j}\mid \mathfrak{d}_i,\mathfrak{e}_i}\mu(\mathfrak{s}_{i,j})$ in the main term yields the main term as equal to
		\begin{align*}
			\frac{\mod{P(a_m,N)}}{\phi({\mathfrak{a}_m\mathfrak{w}})}\sum_{\substack{\mathfrak{u}_1,\ldots,\mathfrak{u}_k\\ \mathfrak{u}_m=1}}\left(\prod_{i=1}^{k}g(\mathfrak{u}_i)\right)\sideset{}{^*}\sum_{\mathfrak{s}_{1,2},\ldots,\mathfrak{s}_{k,k-1}}\left(\prod_{\substack{1\leq i,j\leq k\\ i\neq j}}\mu(\mathfrak{s}_{i,j})\right)\sum_{\substack{\mathfrak{d}_1,\ldots,\mathfrak{d}_k\\\mathfrak{e}_1,\ldots,\mathfrak{e}_k\\\mathfrak{u}_i\mid \mathfrak{d}_i,\mathfrak{e}_i\forall i\\\mathfrak{s}_{i,j}\mid \mathfrak{d}_i,\mathfrak{e}_i\forall i\neq j\\
					\mathfrak{d}_m=\mathfrak{e}_m=1}}\frac{\lambda_{\mathfrak{d}_1,\ldots,\mathfrak{d}_k}\lambda_{\mathfrak{e}_1,\ldots,\mathfrak{e}_k}}{\prod_{i=1}^{k}\phi(\mathfrak{d}_i)\phi(\mathfrak{e}_i)},
		\end{align*}
		where $*$ in the sum means that $\mathfrak{s}_{i,j}$ is coprime to each of $\mf{u}_i, \mf{u}_j, \mathfrak{s}_{a,j}, \mathfrak{s}_{i,b}$ for all $i\neq a, j\neq b$. By the change of variables \eqref{def_ym}, the above main term becomes 
		\begin{align*}
			\frac{\mod{P(a_m,N)}}{\phi(\mathfrak{a}_m\mathfrak{w})}\sum_{\substack{\mathfrak{u}_1,\ldots,\mathfrak{u}_k\\\mathfrak{u}_m=1}}\left(\prod_{i=1}^{k}\frac{\mu(\mathfrak{u}_i)^2}{g(\mathfrak{u}_i)}\right)\sideset{}{^*}\sum_{\mathfrak{s}_{1,2},\ldots,\mathfrak{s}_{k,k-1}}\left(\prod_{\substack{1\leq i,j\leq k\\ i\neq j}}\frac{\mu(\mathfrak{s}_{i,j})}{g(\mathfrak{s}_{i,j})^2}\right)y^{(m)}_{\mathfrak{b}_1,\ldots,\mathfrak{b}_k}y^{(m)}_{\mathfrak{c}_1,\ldots,\mathfrak{c}_k},
		\end{align*}
		where $\mathfrak{b}_i:=\mathfrak{u}_i\prod_{j\neq i}\mathfrak{s}_{i,j}$ and $ \mathfrak{c}_j:=\mathfrak{u}_j\prod_{i\neq j}\mathfrak{s}_{i,j}$. As in Lemma~\ref{Key_lemma1}, we see that there is no contribution to the above sum for $(\mathfrak{s}_{i,j},\mathfrak{w})\neq 1$, and we can take $(\mathfrak{s}_{i,j},\mathfrak{w})= 1$. So we have only two cases: either $\mathfrak{s}_{i,j}=1, \forall i\neq j$, or $\mod{\mathfrak{s}_{i,j}}>D_0$ for some $i\neq j$. For $\mod{\mathfrak{s}_{i,j}}>D_0$, contribution of the above main term is 
		\begin{align*}
			&\ll \frac{(y^{(m)}_{\max})^2\mod{A(N)}}{\phi(\mathfrak{w})\log N}\left(\sum_{\substack{\mod{\mathfrak{u}}<R\\(\mathfrak{u},\mathfrak{w})=1}}\frac{\mu(\mathfrak{u})^2}{g(\mathfrak{u})}\right)^{k-1} \left(\sum_{\mod{\mathfrak{s}_{i,j}}>D_0}\frac{\mu(\mathfrak{s}_{i,j})^2}{g(\mathfrak{s}_{i,j})^2}\right)\left(\sum_{\mathfrak{s}\subseteq A}\frac{\mu(\mathfrak{s})^2}{g(\mathfrak{s})^2}\right)^{k^2 -k-1}\\
			&\ll \frac{(y^{(m)}_{\max})^2\mod{A(N)}\phi(\mathfrak{w})^{k-2} (\log R)^{k-1}}{\mod{\mathfrak{w}}^{k-1}D_0\log N} \ll \frac{(y^{(m)}_{\max})^2\phi(\mathfrak{w})^{k-2} \mod{A(N)}(\log N)^{k-2}}{\mod{\mathfrak{w}}^{k-1}D_0},
		\end{align*}
		where the estimate $\frac{\mod{P(a_m,N)}}{\phi(\mathfrak{a}_m)}\ll \frac{\mod{A(N)}}{\log N}$ follows from the fact that $a_m$ is fixed with respect to $N$ (also, it follows from this Case 1 and Case 2 below). Therefore, we have 
		\begin{align}\label{final_S^m_2}
			S^{(m)}_2=\frac{\mod{P(a_m,N)}}{\phi(\mathfrak{a}_m\mathfrak{w})}\sum_{\mathfrak{u}_1,\ldots,\mathfrak{u}_k}\frac{(y^{(m)}_{\mathfrak{u}_1,\ldots,\mathfrak{u}_k})^2}{\prod_{i=1}^{k}g(\mathfrak{u}_i)}&+O\left(\frac{(y^{(m)}_{\max})^2\phi(\mathfrak{w})^{k-2} \mod{A(N)}(\log N)^{k-2}}{\mod{\mathfrak{w}}^{k-1}D_0}\right)\nonumber\\&+O\left(\frac{y^2_{\max}\mod{A(N)}}{(\log N)^B}\right).
		\end{align}
		Next we evaluate the reduced main term by considering the  function field case and the number field case separately.\\}
		\textbf{Case 1:} \textnormal{Let $A=\mathbb{F}_q[t]$ and degree of $a_m$, deg$(a_m)=D_m$. Then, for $N=q^n$, we have $\mod{A(N)}=N$ and also $P(a_m,N)= P(\mod{\mathfrak{a}_m}N)$. Therefore, it follows from Lemma~\ref{PNT_functionfield}, 
		\begin{align*}
			\mod{P(\mod{\mathfrak{a}_m}N)}&=\frac{q^{D_m} N}{\log_{q}(q^{D_m} N)}+O\left(\frac{(q^{D_m}N)^{1/2}}{\log_{q}(q^{D_m}N)}\right)\\
			&=\frac{q^{D_m} N}{\log_{q}N}\left(1+O\left(\frac{1}{\log_{q}N}\right)\right)+O\left(\frac{N^{1/2}}{\log_{q}N}\right)\\
			&=\mod{\mathfrak{a}_m}\mod{P(N)}+O\left(\frac{N}{(\log_{q}N)^2}\right).
		\end{align*}
 Note that, for sufficiently large $N$, each prime factor of $\mf{a}_m$ is also a prime factor of $\mf{w}$ and in this case $\phi(\mf{a}_m\mf{w})=|\mf{a}_m|\phi(\mf{w})$. So
		\begin{align*}
			\frac{\mod{P(\mod{\mathfrak{a}_m}N)}}{\phi(\mathfrak{a}_m\mf{w})}&=\frac{\mod{P(N)}}{\phi(\mathfrak{w})}+O\left(\frac{N}{\phi(\mf{w})(\log_{q}N)^2}\right).
		\end{align*}
		Also
		\begin{align*}
			&\frac{N}{\phi(\mathfrak{w})(\log_{q}N)^2}\sum_{\mathfrak{u}_1,\ldots,\mathfrak{u}_k}\frac{(y^{(m)}_{\mathfrak{u}_1,\ldots,\mathfrak{u}_k})^2}{\prod_{i=1}^{k}g(\mathfrak{u}_i)}\ll \frac{N(y^{(m)}_{\max})^2}{\phi(\mathfrak{w})(\log_{q}N)^2}\left(\sum_{\mod{\mathfrak{u}}<R}\frac{\mu(\mathfrak{u})^2}{g(\mathfrak{u})}\right)^{k-1}\\
			&\ll \frac{(y^{(m)}_{\max})^2\mod{A(N)}\phi(\mathfrak{w})^{k-2} (\log N)^{k-3}}{\mod{\mathfrak{w}}^{k-1}D_0}.
		\end{align*}}
		\textbf{Case 2:} \textnormal{Let $A=\mathcal{O}_K$ with $[K:\mathbb{Q}]=d$ and signature $(r_1, r_2)$. Then by Lemma~\ref{GPNT} for $A_0(\mathfrak{a}_m,N)$ (see \eqref{A_0(am,N)}), we have the number of primes in $A_0(\mathfrak{a}_m,N)$,
		\begin{align*}
\mod{P_0(a_m,N)}=&m_K\int_{2}^{|\sigma_1(a_m)|N}\hspace{-18pt}\ldots\int_{2}^{|\sigma_{r_1}(a_m)|N}\hspace{-8pt}\int_{2}^{(|\sigma_{r_1+1}(a_m)|N)^2}\hspace{-18pt}\ldots\int_{2}^{(|\sigma_{r_1+r_2}(a_m)|N)^2}\hspace{-11pt}\frac{du_1\ldots du_{r_1+r_2}}{\log(u_1\ldots u_{r_1+r_2})}\\
			&+O\left(|\mathfrak{a}_m|N^d e^{-c\sqrt{\log(|\mathfrak{a}_m| N^d)}}\right),
		\end{align*}
		where $\sigma_j$ are real embeddings for $1\leq j\leq r_1$ and $\sigma_j$ are complex embeddings for $(r_1+1)\leq j\leq (r_1+r_2)$. Then by Lemma~\ref{Value_of_int}, 
		\begin{align*}
			\mod{P_0(a_m,N)}&= \frac{ m_K|\mathfrak{a}_m|N^d}{\log(|\mathfrak{a}_m|N^d)} +O\left(\frac{ m_K|\mathfrak{a}_m|N^d}{(\log(|\mathfrak{a}_m|N^d))^2}\right)+O\left(|\mathfrak{a}_m|N^d e^{-c\sqrt{\log(|\mathfrak{a}_m| N^d)}}\right)\\
			&= \frac{ m_K|\mathfrak{a}_m|N^d}{\log(|\mathfrak{a}_m|N^d)}+O\left(\frac{N^d}{(\log N)^2}\right)\\
			&=\frac{ m_K|\mathfrak{a}_m|N^d}{\log(N^d)}\left(1+O\left(\frac{1}{\log N}\right)\right)+O\left(\frac{N^d}{(\log N)^2}\right)\\
			&=|\mathfrak{a}_m||P_0(N)|+O\left(\frac{N^d}{(\log N)^2}\right).
		\end{align*}
		Since $P(N)=P_0(2N)\setminus P_0(N)$, we have
		\begin{align*}
			|P(a_m,N)|=|\mathfrak{a}_m||P(N)|+O\left(\frac{N^d}{(\log N)^2}\right)
		\end{align*}
		and by similar argument as in Case~1, 
		\begin{align*}
			\frac{|P(a_m,N)|}{\phi(\mathfrak{a}_m\mf{w})}&=\frac{|P(N)|}{\phi(\mathfrak{w})}+O\left(\frac{N^d}{\phi(\mf{w})(\log N)^2}\right).
		\end{align*}
		Also, as we know $N^d=O(|A(N)|)$ (for example, see proof of Corollary 2.8 in \cite{Castillo}), it follows that 
		\begin{align*}
			&\frac{N^d}{\phi(\mathfrak{w})(\log N)^2}\sum_{\mathfrak{u}_1,\ldots,\mathfrak{u}_k}\frac{(y^{(m)}_{\mathfrak{u}_1,\ldots,\mathfrak{u}_k})^2}{\prod_{i=1}^{k}g(\mathfrak{u}_i)}\\
			&\ll \frac{|A(N)|(y^{(m)}_{\max})^2}{\phi(\mathfrak{w})(\log N)^2}\left(\sum_{\mod{\mathfrak{u}}<R}\frac{\mu(\mathfrak{u})^2}{g(\mathfrak{u})}\right)^{k-1}\\
			&\ll \frac{(y^{(m)}_{\max})^2\mod{A(N)}\phi(\mathfrak{w})^{k-2} (\log N)^{k-3}}{\mod{\mathfrak{w}}^{k-1}D_0}.
		\end{align*}}
		
		\noindent \textnormal{From Case~1, Case~2 and from \eqref{final_S^m_2}, we obtain \eqref{S_2m}, which completes the proof.} \qed
	\end{proof2}

	To prove Propositions~\ref{1st_prop} and Propositions~\ref{2nd_prop}, we need the following lemmas.
	\begin{lemma}[Lemma 2.6 in \cite{Castillo}]
		If $\mathfrak{r}_m=1$, the trivial ideal, then
		\begin{align*}
			y^{(m)}_{\mathfrak{r}_1,\ldots,\mathfrak{r}_k}=\sum_{\mathfrak{t}_m}\frac{y_{\mathfrak{r}_1,\ldots,\mathfrak{r}_{m-1},\mathfrak{t}_m,\mathfrak{r}_{m+1},\ldots,\mathfrak{r}_k}}{\phi(\mathfrak{t}_m)}+O\left(\frac{y_{\max}\phi(\mathfrak{w})\log R}{|\mathfrak{w}|D_0}\right).	
		\end{align*}
	\end{lemma}

	We now choose\footnote{See  Section~6 in \cite{Maynard} and  Section~2 in \cite{Castillo}.} 
	\begin{align}\label{smooth_y}
		y_{\mathfrak{r}_1,\ldots,\mathfrak{r}_k}:=F\left(\frac{\log |\mathfrak{r}_1|}{\log R},\ldots,\frac{\log |\mathfrak{r}_k|}{\log R}\right),
	\end{align}
	where $F:[0,1]^k\to \mathbb{R}$ is a piecewise differentiable function supported on $\mathcal{R}_k=\{(x_1,\ldots,x_k)\in[0,1]^k:\sum_{i=1}^{k}x_i\leq 1\}.$ Also, we let $y_{\mathfrak{r}_1,\ldots,\mathfrak{r}_k}$ to be zero unless $ \mathfrak{r}=\prod_{i=1}^{k}\mathfrak{r}_i$  is squarefree and coprime to $\mathfrak{w}$.
	\begin{lemma}[Lemma 2.7 in \cite{Castillo}]\label{Sum_to_integral}
		Let $\gamma$ be a multiplicative function on nonzero ideals of $A$ such that there are constants $\kappa>0, A_1>0, A_2\geq 1,$ and $L\geq 1$  satisfying 
		\begin{align*}
			0\leq \frac{\gamma(\mathfrak{p})}{|\mathfrak{p}|}\leq 1-A_1,
		\end{align*}
		and 
		\begin{align*}
			-L\leq \sum_{w<|\mathfrak{p}|\leq z}\frac{\gamma(\mathfrak{p})\log |\mathfrak{p}|}{|\mathfrak{p}|}-\kappa\log\frac{z}{w}\leq A_2
		\end{align*}
		for any $2\leq w\leq z$. Let $g$ be totally multiplicative  with $g(\mathfrak{p})=\frac{\gamma(\mathfrak{p})}{|\mathfrak{p}|-\gamma(\mathfrak{p})}$ and $G:[0,1]\to\mathbb{R}$ be piecewise differentiable, and $G_{\max}=\sup_{t\in[0,1]}(|G(t)|+|G'(t)|).$ Then 
		\begin{align*}
			\sum_{|\mathfrak{d}|<z}\mu(\mathfrak{d})^2g(\mathfrak{d})G\left(\frac{\log |\mathfrak{d}|}{\log z}\right)=\mathfrak{S}\frac{c^{\kappa}_A(\log z)^\kappa}{\Gamma(\kappa)}\int_{0}^{1}G(x)x^{\kappa -1}dx+O_{A,A_1,A_2,\kappa}(\mathfrak{S}LG_{\max}),
		\end{align*}
		where $c_A:= \text{Res}_{s=1}\zeta_A(s)$ for $\zeta_A(s)$ in \eqref{zeta_function} and
		\begin{align*}
			\mathfrak{S}=\prod_{\mathfrak{p}}\left(1-\frac{\gamma(\mathfrak{p})}{|\mathfrak{p}|}\right)^{-1}\left(1-\frac{1}{|\mathfrak{p}|}\right)^\kappa.
		\end{align*}	
	\end{lemma}
	By proceeding in a similar manner as in the proofs of Lemma~6.2 and 6.3 in \cite{Maynard} and with repeated application of the above lemma with $\kappa=1$, we obtain the following lemma.
	\begin{lemma}\label{last_lemma}
		Let $y_{\mathfrak{r}_1,\ldots,\mathfrak{r}_k}$ be given in \eqref{smooth_y} and 
		\begin{align*}
			F_{\max}:=\sup_{(t_1,\ldots,t_k)\in[0,1]^k}|F(t_1,\ldots,t_k)|+\sum_{i=1}^{k}|\frac{\partial F}{\partial t_i}(t_1,\ldots,t_k)|
		\end{align*}
		Then, we have 
		\begin{align*}
			S_1=\frac{\phi(\mathfrak{w})^k |A(N)| (c_A\log R)^k}{|\mathfrak{w}|^{k+1}}I_k(F)+O\left(\frac{F^2_{\max}\phi(\mathfrak{w})^k |A(N)|(\log R)^k}{|\mathfrak{w}|^{k+1}D_0}\right),
		\end{align*}
		and 
		\begin{align*}
			S^{(m)}_2= \frac{\phi(\mathfrak{w})^k |P(N)|(c_A\log R)^{k+1}}{|\mathfrak{w}|^{k+1}}J^{(m)}_k(F)+O\left(\frac{F^2_{\max}\phi(\mathfrak{w})^k |A(N)|(\log R)^k}{|\mathfrak{w}|^{k+1}D_0}\right),
		\end{align*}
		where $I_k(F)$ and $J_k(F)$ are as in \eqref{I_k} and \eqref{J_k}, respectively and $c_A$ as in Lemma~\ref{Sum_to_integral}. 
	\end{lemma}
	\section{Proof of the Results}
	We note that Proposition~\ref{1st_prop} follows immediately from  Lemma~\ref{last_lemma}. Below we prove Proposition~\ref{2nd_prop}.\\
	
	\noindent\emph{Proof of Proposition~\ref{2nd_prop}.} Let us recall that $R=|A(N)|^{\theta/2 -\delta}$, for some fixed small $\delta>0$ and 
	\begin{align*}
		M_k=\sup_{F\in\mathcal{S}_k}\frac{\sum_{m=1}^{k}J^{(m)}_k(F)}{I_k(F)},
	\end{align*}
	where $\mathcal{S}_k$ be set of all piecewise differentiable functions $F:[0,1]^k\to\mathbb{R}$ supported on the simplex $\mathcal{R}_k=\{(t_1,\ldots,t_k)\in[0,1]^k:\sum_{i=1}^{k} t_i\leq1\}$ with $I_k(F)\neq 0$  and $J^{(m)}_k(F)\neq 0$, for each $m$. Choose $F_0\in \mathcal{S}_k$ such that $\sum_{m=1}^{k}J^{(m)}_k(F_0)>(M_k-\delta)I_k(F_0)$. Then  choosing sieve weights $\lambda$ in terms of $F_0$ and using Lemma~\ref{last_lemma}, we have for any positive real number $\rho>0$,
	\begin{align*}
		S:&=S_2-\rho S_1\\
		&= \frac{\phi(\mathfrak{w})^k|A(N)|(c_A\log R)^k}{|\mathfrak{w}|^{k+1}}\left(\frac{(c_A\log R) |P(N)|}{|A(N)|}\sum_{m=1}^{k}J^{(m)}_k(F_0)-\rho I_k(F_0)+o(1)\right)\\
		&\geq \frac{\phi(\mathfrak{w})^k|A(N)|(c_A\log R)^k}{|\mathfrak{w}|^{k+1}}I_k(F_0)\left(\Delta\cdot\left(\frac{\theta}{2}-\delta\right)\left(M_k-\delta\right)-\rho+o(1)\right),
	\end{align*}
	where 
	\begin{align*}
		\Delta:=c_A\cdot\lim_{N\to\infty}\frac{|P(N)|\log|A(N)|}{|A(N)|}.
	\end{align*}
	From the proof of Corollary~2.8 in \cite{Castillo}, we know that the above limit exists and $\Delta=1$. Thus,
	\begin{align*}
		S&\geq \frac{\phi(\mathfrak{w})^k|A(N)|(c_A\log R)^k}{|\mathfrak{w}|^{k+1}}I_k(F_0)\left(\left(\frac{\theta}{2}-\delta\right)\left(M_k-\delta\right)-\rho+o(1)\right)\\
		&=\frac{\phi(\mathfrak{w})^k|A(N)|(c_A\log R)^k}{|\mathfrak{w}|^{k+1}}I_k(F_0)\left(\frac{\theta M_k}{2}-\rho -\delta(M_k+\theta-\delta)+o(1)\right).
	\end{align*}
	If we choose $\rho =\frac{\theta M_k}{2}-\epsilon$, then for a suitable choice of $\delta$ depending on $\epsilon$, we have $S>0$ and this implies that there are infinitely many $\alpha\in A$ such that at least $\lfloor\rho+1\rfloor$ of the $a_i\alpha+h_i, 1\leq i\leq k$ are prime. Now for suitably chosen $\epsilon,$ we have 
 \begin{equation*}
     \lfloor\rho+1\rfloor= \left\lceil\frac{\theta M_k}{2}\right\rceil.
 \end{equation*} This completes the proof.\qed\\
	
\noindent\emph{Proof of the Theorem~\ref{Thm1} and Theorem~\ref{Thm2}.} From Lemma~\ref{M_k}, we have for sufficiently large $k$, $M_k>\log k-2\log\log k-2$ and the later quantity is $>4m'$, if $k\geq C(m')^2 e^{4m'}$ for some absolute constant $C$ (see the proof of Proposition 4.3 in \cite{Maynard}). Therefore,
	\begin{align*}
		\frac{\theta M_k}{2}> 2\theta m'\geq m,
	\end{align*}
	if $k\geq C(\frac{m}{2\theta})^2e^{\frac{2m}{\theta}}$. This completes the proofs with $k_0=\lceil C(\frac{m}{2\theta})^2e^{\frac{2m}{\theta}}\rceil$.\qed
	
	\begin{remark}
		\textnormal{(i) As for number field $K$, the level of distribution $\theta$ of primes depend on the number field $K$ and as for the function field $\mathbb{F}_q(t)$, $\theta$ does not depend on $\mathbb{F}_q(t)$, we have $k_0=k_0(m,K)$ and $k_0=k_0(m)$, respectively.\\
		(ii) We know by Theorems~\ref{Hinz} and \ref{Hayes}, for totally real number fields $K$ and for any function field $\mathbb{F}_q(t)$, the level of distribution is $\theta $ for any $\theta<1/2$ . So taking $\theta=\frac{1}{2}-\frac{1}{k}$, we see that}
		\begin{align*}
			\frac{\theta M_k}{2}&>\left(\frac{1}{4}-\frac{1}{2k}\right)(\log k-2\log\log k-2)= \left(1/4 +o_{k\to\infty}(1)\right)\log k.
		\end{align*}
		\textnormal{Therefore, for these fields, for sufficiently large $k$, any admissible set $\mathcal{H}=\{a_1\alpha+h_1,\ldots,a_k\alpha+h_k\}$ contains at least $\left(1/4 +o_{k\to\infty}(1)\right)\log k$ prime components for infinitely many $\alpha\in A$, which is also mentioned by Maynard \cite{Maynard}, for $A=\mathbb{Z}$.}
	\end{remark}
\noindent\emph{Proof of Corollary~\ref{main_cor}.}	Let $m\geq 2$ and $k\geq k_0(m,K)$, where $k_0(m,K)$ is as in Theorem~\ref{Thm1}. Let $\{\alpha+h_1, \alpha+h_2,\ldots, \alpha+h_k\}$ be an admissible set in $\mc{O}_K$. Let $b_i=hh_i+b$ for $1\leq i\leq k$ and let 
\begin{align*}
	\mc{H}_1:=\{h\alpha+b_1, h\alpha+b_2,\ldots, h\alpha+b_k\}.
\end{align*} 
We claim that $\mc{H}_1$ is also an admissible set. Cleary, for any prime ideal $\mf{p}\mid h$, 
\begin{align*}
	\mf{p}\nmid \prod_{i=1}^{k}(h\alpha+b_i):=L(\alpha),
\end{align*} 
for any $\alpha\in\mc{O}_K$, since $h$ and $b$ are coprime. So let $\mf{p}$ be any prime ideal in $\mc{O}_K$ and $\mf{p}\nmid h$. Since $\{\alpha+h_1, \alpha+h_2,\ldots, \alpha+h_k\}$ is admissible, there exists $\alpha_\mf{p}\in\mc{O}_K$ such that 
\begin{align}\label{1st_eq}
\mf{p}\nmid \prod_{i=1}^{k}(\alpha_\mf{p}+h_i).
\end{align}
Suppose 
\begin{align}\label{y}
	\mf{p}\mid \prod_{i=1}^{k}(h\alpha+b_i) 
\end{align}
for all $\alpha\in\mc{O}_K.$ We choose $\alpha=-\tilde{h}b+\alpha_\mf{p}$ where $h\tilde{h}\equiv 1(\bmod \mf{p})$. Inserting this $\alpha$ in \eqref{y}, we have 
\begin{align*}
	\mf{p}\mid \prod_{i=1}^{k}(\alpha_\mf{p}+h_i),
\end{align*}
which contradicts \eqref{1st_eq}. 
Therefore, there exists at least one $\beta_\mf{p}\in\mc{O}_K$ such that 
\begin{align*}
    \mf{p}\nmid \prod_{i=1}^{k}(h\beta_\mf{p}+b_i) 
\end{align*}
and hence the set $\mc{H}_1$ is admissible.
Then applying  Theorem~\ref{Thm1} to the admissible set $\mc{H}_1$ for given $m$, we obtain that there are infinite many $\alpha\in\mc{O}_K$ such that at least $m$ of $h\alpha+hh_i+b, 1\leq i\leq k$ are prime. Also $h\alpha+hh_i+b$ lies in the congruence class $b (\bmod \mf{h}).$ Therefore there are infinitely many $r\in\mb{N}$ such that 
\begin{align*}
	\gamma_{r1}\equiv\gamma_{r2}\equiv\cdots\equiv\gamma_{rm}\equiv b (\bmod \mathfrak{h}),
\end{align*}
where $\gamma_{rj}$ denote primes in $\mc{O}_K$.\qed\\

\noindent\emph{Proof of Corollary~\ref{main_cor_2}}. The proof follows in a similar manner to Corollary~\ref{main_cor}
by Theorem~\ref{Thm2}. \qed



\end{document}